\documentclass[12pt, a4paper]{article}
\usepackage[dvips]{graphicx}
\usepackage[T1]{fontenc}
\usepackage[latin1]{inputenc}
\usepackage[english]{babel}
\usepackage{times}
\usepackage{amsmath}
\usepackage{amssymb}
\usepackage{amscd}

\usepackage{latexsym}
\usepackage{tabularx}
\usepackage{here}
\usepackage{color}
\usepackage{vmargin}
\setmarginsrb{2cm}{2cm}{2cm}{2cm}{0cm}{2cm}{0cm}{1cm}
\begin{document}

\title{\Large{\bf $(3,1)^*$-choosability of planar graphs without adjacent short cycles \\}}

\author{Min Chen$^{a}$\thanks{Research supported by NSFC (No.11101377). Email:
chenmin@zjnu.cn}, \ Andr\'{e} Raspaud$^{b}$\thanks{ Research partially
supported by ANR-NSC Project GRATEL - ANR-09-blan-0373-01 and
NSC99-2923-M-110-001-MY3. Email:
andre.raspaud@labri.fr. Tel: +33 5 40 0 69
29. Fax: +33 5 40 00 66 69.}\\
\small \ \ $^{a}$ \ {\em Department of Mathematics, Zhejiang Normal
University, Jinhua 321004, China}\\ \
\small \ \ $^{b}$ \ {\em  LaBRI UMR CNRS 5800, Universite Bordeaux
I, \small33405 Talence Cedex, France. 
}}

\maketitle
\newcommand{\ch}{\omega}
\newcommand{\nch}{\omega^*}
\newcommand{\qed}{\hfill $\Box$ }
\newtheorem{corollary}{Corollary}
\newtheorem{definition}{Definition}
\newtheorem{question}{Question}
\newtheorem{proposition}{Proposition}
\newtheorem{theorem}{Theorem}
\newtheorem{lemma}{Lemma}
\newtheorem{conjecture}{Conjecture}
\newtheorem{sketch}{Sketch of proof}
\newtheorem{ndm}{N.d.M.}
\newtheorem{observation}{Observation}
\newtheorem{idea}{Idea}
\newtheorem{remark}{Remark}
\newtheorem{claim}{Claim}
\newtheorem{example}{Example}

\newcommand{\proof}{\noindent{\bf Proof.}\ \ }
\baselineskip=17pt
\parindent=0.5cm

\begin{abstract}

A list assignment of a graph $G$ is a function $L$ that assigns a list $L(v)$
of colors to each vertex $v\in V(G)$.
An $(L,d)^*$-coloring is a mapping $\pi$ that assigns a color $\pi(v)\in L(v)$ to each
vertex $v\in V(G)$ so that at most $d$ neighbors of $v$ receive color $\pi(v)$.
A graph $G$ is said to be $(k,d)^*$-choosable
if it admits an $(L,d)^*$-coloring for every list assignment $L$ with $|L(v)|\ge k$ for all $v\in V(G)$.
In 2001, Lih et al. \cite{LSWZ-01} proved that planar graphs without $4$- and $l$-cycles
are $(3,1)^*$-choosable, where $l\in \{5,6,7\}$. Later, Dong and Xu \cite{DX-09} proved that
planar graphs without $4$- and $l$-cycles are $(3,1)^*$-choosable,
where $l\in \{8,9\}$.

There exist planar graphs containing $4$-cycles that are not $(3,1)^*$-choosable
(Crown, Crown and Woodall, 1986 \cite{CCW-86}).
This partly explains the fact that in all above known sufficient
conditions for the $(3,1)^*$-choosability of planar graphs the $4$-cycles are completely
forbidden.
In this paper we allow $4$-cycles nonadjacent to relatively short cycles. More precisely,
we prove that every planar graph without $4$-cycles adjacent to $3$- and $4$-cycles
is $(3,1)^*$-choosable. This is a common strengthening of all above
mentioned results.
Moreover as a consequence we give a  partial answer to a question of Xu and Zhang
\cite{XZ-07} and show that every planar graph without $4$-cycles is
$(3,1)^*$-choosable.

\bigskip

\medskip
\end{abstract}

\section{Introduction}\label{intro}

\baselineskip=20pt

All graphs considered in this paper are finite, loopless, and
without multiple edges. A {\em plane graph} is a particular drawing
of a planar graph in the Euclidean plane. For a graph $G$, we use
$V(G)$, $E(G)$, $|G|$, $|E(G)|$ and $\delta(G)$ to denote its
vertex set, edge set, order, size and minimum degree, respectively.
For $v\in V(G)$, $N_G(v)$ denotes the set of neighbors of $v$ in
$G$. If there is no confusion about the context, we write $N(v)$ for
$N_G(v)$.

A $k$-coloring of $G$ is a mapping $\pi$ from $V(G)$ to a color set $\{1,2,\cdots, k\}$
such that $\pi(x)\neq \pi(y)$ for any adjacent vertices $x$ and $y$.
A graph is $k$-colorable if it has a $k$-coloring.
Cowen, Cowen, and Woodall \cite{CCW-86} considered {\em defective} colorings of
graphs. A graph $G$ is said to be
{\em $d$-improper $k$-colorable}, or simply, {\em $(k,d)^*$-colorable},
if the vertices of $G$ can be colored with $k$ colors in such a way that each vertex has at most $d$ neighbors
receiving the same color as itself.
Obviously, a $(k,0)^*$-coloring is an ordinary proper $k$-coloring.

A {\em list assignment} of $G$ is a function $L$ that assigns a list $L(v)$
of colors to each vertex $v\in V(G)$.
An $L$-coloring with
impropriety of integer $d$, or simply an {\em $(L,d)^*$-coloring}, of $G$ is a mapping
$\pi$ that assigns a color $\pi(v)\in L(v)$ to each
vertex $v\in V(G)$ so that at most $d$ neighbors of $v$ receive color $\pi(v)$.
A graph is {\em $k$-choosable} with impropriety of integer $d$,
or simply {\em $(k,d)^*$-choosable}, if there exists an $(L,d)^*$-coloring for every list assignment $L$ with $|L(v)|\ge k$
for all $v\in V(G)$. Clearly, a $(k,0)^*$-choosable is the ordinary $k$-choosability introduced by
Erd\H{o}s, Rubin and Taylor \cite{ERT-79}
and independently by Vizing \cite{Vizing-76}.

The concept of list improper coloring was independently introduced by \v{S}krekovski \cite{Srekovski-99}
and Eaton and Hull \cite{EH-99}. They proved that every planar graph is $(3,2)^*$-choosable and
every outerplanar graph is $(2,2)^*$-choosable.
These are both improvement of the results showed in \cite{CCW-86}
which say that
every planar graph is $(3,2)^*$-colorable and
every outerplanar graph is $(2,2)^*$-colorable.
Let $g(G)$ denote the {\em girth} of a graph $G$, i.e.,
the length of a shortest cycle in $G$.
The $(k,d)^*$-choosability of planar graph $G$ with given $g(G)$ has
been studied by \v{S}krekovski in \cite{Srekovski-00}. He proved that
every planar graph $G$ is $(2,1)^*$-choosable if $g(G)\ge 9$, $(2,2)^*$-choosable if $g(G)\ge 7$,
$(2,3)^*$-choosable if $g(G)\ge 6$, and $(2,d)^*$-choosable if $d\ge 4$ and $g(G)\ge 5$.
Recently, Cushing and Kierstead \cite{CK-10} proved that
every planar graph is $(4,1)^*$-choosable.
So it would be interesting to investigate the sufficient
conditions of $(3,1)^*$-choosability of subfamilies of planar graphs where some families of cycles are forbidden.
\v{S}krekovski proved in \cite{Srekovski-99-1}
that every planar graph without $3$-cycles is $(3,1)^*$-choosable.
Lih et al. \cite{LSWZ-01} proved that planar graphs without $4$- and $l$-cycles
are $(3,1)^*$-choosable, where $l\in \{5,6,7\}$. Later, Dong and Xu \cite{DX-09} proved that
planar graphs without $4$- and $l$-cycles are  $(3,1)^*$-choosable, where $l\in \{8,9\}$.
Moreover,  Xu and Zhang \cite{XZ-07} asked the following question:

\begin{question}\label{Q-1}
Is it true that every planar graph without adjacent triangles is $(3,1)^*$-choosable?
\end{question}

Recall that there is a planar graph containing $4$-cycles that is
not $(3,1)^*$-colorable \cite{CCW-86}.
Therefore, while describing
$(3,1)^*$-choosability planar graphs,
one must impose these or those restrictions on $4$-cycles. Note that in all previously
known sufficient conditions for the $(3,1)^*$-choosability of planar graphs,
the $4$-cycles are completely forbidden. In this
paper we allow $4$-cycles, but disallow them to have a common edge with relatively short cycles.

The purpose of this paper is to prove the following

\begin{theorem}\label{main01}
Every planar graph without $4$-cycles adjacent to $3$- and $4$-cycles
is $(3,1)^*$-choosable.
\end{theorem}

Clearly, Theorem \ref{main01} implies Corollary \ref{cor1} which
is a common strengthening of the results in \cite{LSWZ-01, DX-09}.

\begin{corollary}\label{cor1}
Every planar graph without $4$-cycles is $(3,1)^*$-choosable.
\end{corollary}

Moreover, Theorem \ref{main01} partially answers Question \ref{Q-1}, since
adjacent triangles can be regarded as a $4$-cycle adjacent to a $3$-cycle.

\section{Notation}

A vertex of degree $k$ (resp. at least $k$, at most $k$) will be
called a {\em $k$-vertex} (resp. {\em $k^+$-vertex}, {\em $k^-$-vertex}). A similar
notation will be used for cycles and faces.
A {\em triangle} is synonymous with a 3-cycle.
For $f\in F(G)$, we use $b(f)$ to denote the boundary walk of $f$ and write
$f=[u_1u_2\cdots u_n]$ if $u_1,u_2,\cdots,u_n$ are the boundary
vertices of $f$ in cyclic  order.
For any $v\in V(G)$, we let $v_1,v_2,\cdots, v_{d(v)}$ denote the
neighbors of $v$ in a cyclic order. Let $f_i$ be the
face with $vv_i$ and $vv_{i+1}$ as two boundary edges for
$i=1,2,\cdots,d(v)$, where indices are taken modulo $d(v)$.
Moreover, we let $t(v)$ denote the number
of $3$-faces incident to $v$ and let $n_3(v)$ denote the
number of $3$-vertices adjacent to $v$.

An $m$-face $f=[v_1v_2\cdots v_m]$ is called an $(a_1,a_2,\cdots, a_m)$-{\em face} if the
degree of the vertex $v_i$ is $a_i$ for $i=1,2,\cdots, m$.
Suppose $v$ is a $4$-vertex incident to a $4^-$-face $f$ and adjacent to two
$3$-vertices not on $b(f)$. If $d(f)=3$, then  we call $v$ a {\em light} $4$-vertex.
Otherwise, we call $v$ a {\em soft} $4$-vertex if $d(f)=4$.
A vertex $v$ is called an {\em $\mathcal{S}$-vertex} if it is either a $3$-vertex or a light $4$-vertex.
Moreover, we say a $3$-face $f=[v_1v_2v_3]$ is an $(a_1,*,a_3)$-face if $d(v_i)=a_i$ for each $i\in \{1,3\}$ and
$v_2$ is an $\mathcal{S}$-vertex.
Suppose $v$ is a $5$-vertex incident to two 3-faces $f_1=[vv_1v_2]$
and $f_3=[vv_3v_4]$. Let $v_5$ be the neighbour of $v$ not belonging to
the $3$-faces. If $d(v_5)=3$ and $f_1$ is a $(5,*,4)$-face,
then we call $v$ a {\em bad} $5$-vertex.

\begin{figure}[ht]
\centering
\includegraphics[width=4in]{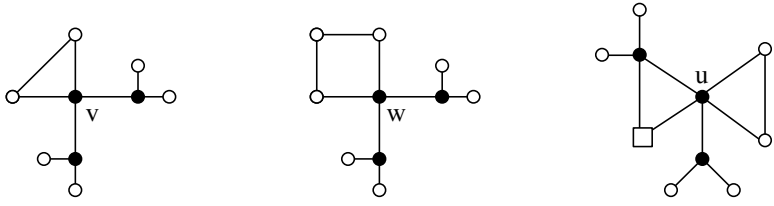}\\
\caption{A light $4$-vertex $v$, a soft $4$-vertex $w$ and a bad $5$-vertex $u$.}
\label{fig2}
\end{figure}

For all figures in the following section, a vertex
is represented by a solid circle when all of its incident edges are
drawn; otherwise it is represented by a hollow circle. Moreover, we use a hollow square
to denote an $\mathcal{S}$-vertex.

\section{Proof of Theorem \ref{main01}}

The proof of Theorem \ref{main01} is done by reducible configurations and discharging procedure.
Suppose the theorem is not true. Let $G$ be a counterexample
with the least number
of vertices and edges embedded in the plane.
Thus, $G$ is connected.
We will apply a discharging procedure to reach
a contradiction.

We first define a weight function $\ch$ on the vertices and faces of $G$
by letting  $\ch(v)=3d(v)-10$ if $v\in V(G)$ and $\ch(f)=2d(f)-10$ if
$f\in F(G)$. It follows from Euler's formula
$|V(G)|-|E(G)|+|F(G)|=2$ and the relation $\sum_{v\in
V(G)}{d(v)}=\sum_{f\in
 F(G)}{d(f)}=2|E(G)|$ that the total sum of weights of the vertices and faces is equal to

\begin{equation*}
\sum\limits_{v\in V(G)}(3d(v)-10)+\sum\limits_{f\in
F(G)}(2d(f)-10)=-20.
\end{equation*}

We then design appropriate discharging rules and redistribute
weights accordingly. Once the discharging is finished, a new weight
function $\nch$ is produced. The total sum of weights is kept fixed
when the discharging is in process. Nevertheless, after the
discharging is complete, the new weight function satisfies
$\nch(x)\geq 0$ for all $x\in V(G)\cup F(G)$. This leads to the
following obvious contradiction,
$$
-20 = \sum_{x\in V(G)\cup F(G)}\ch(x) = \sum_{x\in V(G)\cup
    F(G)}\nch(x) \ge 0
$$
and hence demonstrates that no such counterexample can exist.

\subsection{Reducible configurations of $G$}

In this section, we will establish structural properties of $G$.
More precisely, we prove that some configurations are reducible.
Namely, they cannot appear in $G$ because of the minimality of $G$.
Since $G$ does not contain a $4$-cycle adjacent to an $i$-cycle,
where $i=3,4$, by hypothesis,
the following fact is easy to observe and will be frequently used throughout this paper
without further notice.

\begin{observation}\label{three structures}
$G$ does not contain the following structures:

\noindent{\rm {(a)}} adjacent $3$-cycles;

\noindent{\rm {(b)}} a $4$-cycle adjacent to a $3$-cycle;

\noindent{\rm {(c)}} a $4$-cycle adjacent to a $4$-cycle.
\end{observation}

We first present Lemma \ref{lem-0}, whose proof
was provided in \cite{LSWZ-01}.

\begin{lemma}\label{lem-0}{\rm\cite{LSWZ-01}}

\noindent{\rm (A1)}  $\delta(G) \geq 3$.

\noindent{\rm (A2)}  No two adjacent $3$-vertices.

\noindent{\rm (A3)} There is no $(3,4,4)$-face.

\end{lemma}

Before showing Lemmas \ref{key-0}-\ref{key},
we need to introduce some useful concepts,
which were firstly defined by Zhang in \cite{Zhang-12}.

\begin{definition}
{\em For $S\subseteq V(G)$,
let $G[S]$ denote the
subgraph of $G$ induced by $S$.
We simply write
$G-S=G[V(G)\setminus S]$.
Let $L$ be an arbitrary list assignment of $G$,
and $\pi$ be an $(L,1)^*$-coloring of $G-S$.
For each $v\in S$, let $L_{\pi}(v)=L(v)\setminus \{\pi(u): u\in N_{G-S}(v)\}$,
and we call $L_{\pi}$ an {\em induced assignment} of $G[S]$ from $\pi$.
We also say that $\pi$ can be extended to $G$ if $G[S]$ admits an $(L_{\pi},1)^*$-coloring.}
\end{definition}

\begin{figure}[ht]
\centering
\includegraphics[width=1.4in]{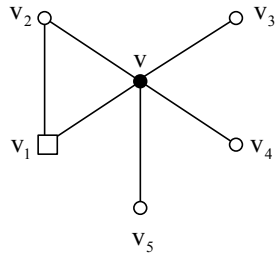}
\caption{The configuration (Q) in Lemma \ref{key-0}.}
\label{fig3}
\end{figure}

\begin{lemma}\label{key-0}
Suppose that $G$ contains the configuration $(Q)$, depicted in Figure \ref{fig3}.
Let $\pi$ be an $(L,1)^*$-coloring of $G-S$, where $S=\{v, v_1, v_2, v_3, v_4\}$.
Denote by $L_{\pi}$ an induced list assignment of $G[S]$.
If $|L_{\pi}(v_i)|\ge 1$ for each $i\in \{1,\cdots, 4\}$,
then $\pi$ can be extended to the whole graph $G$.
\end{lemma}

\proof Since $|L_{\pi}(v_i)|\ge 1$ for each $i\in \{1,\cdots, 4\}$,
we can color each $v_i$ with a color $\pi(v_i)\in L_{\pi}(v_i)$ properly.
Note that $|L_{\pi}(v)|\ge 2$.
If there exists a color in $L_{\pi}(v)$ which appears at most once on the set $\{v_1, v_2,v_3, v_4\}$,
then we assign such a color to $v$. It is easy to check that the resulting
coloring is an $(L,1)^*$-coloring and thus we are done.
Otherwise, w.l.o.g., suppose $L(v)=\{1,2,3\}$, $\pi(v_5)=1$, and each color in $\{2,3\}$
appears exactly twice on the set $\{v_1, v_2, v_3,v_4\}$.
W.l.o.g., suppose $\pi(v_1)=2$.

By definition, we see that $v_1$ is either a $3$-vertex
or a light $4$-vertex. We label two steps in the proof
for future reference.

(i)\ If $d(v_1)=3$, then $|L_{\pi}(v_1)|\ge 2$. We may assign color $2$ to $v$
and then recolor $v_1$ with a color in $L_{\pi}(v_1)\setminus \{2\}$.

(ii)\ If $v_1$ is a light $4$-vertex, denote by
$x_1, y_1$ the other two neighbors which are different from $v$ and $v_2$.
Erase the color of $v_1$, color $v$ with $2$,
and recolor $x_1$ and $y_1$ with a color different from
its neighbors. We can do this since $d(x_1)=d(y_1)=3$ by definition.
Next, we will show how to extend the resulting coloring, denoted by $\pi'$, to $G$.
If $\pi'(v_2)\notin \{\pi'(x_1), \pi'(y_1)\}$,
then color $v_1$ with a color in $L(v_1)\setminus \{2, \pi'(x_1)\}$.
Otherwise, we color $v_1$ with a color in $L(v_1)\setminus \{2, \pi'(v_2)\}$.
In each case, one can easily check that the obtained coloring
of $G$ is an $(L,1)^*$-coloring.

Therefore, we complete the
proof of Lemma \ref{key-0}.\qed

\begin{lemma}\label{lem-1}
$G$ satisfies the following.

\noindent{\rm (B1)}\ A $4$-vertex is adjacent to at most two $3$-vertices.

\noindent{\rm (B2)}\ There is no $(4^-,4^-,4^-)$-face.

\noindent{\rm (B3)}\ There is no $(5^+,4,4)$-face which is incident to two light $4$-vertices.

\noindent{\rm (B4)}\ There is no $5$-vertex incident to a $(5,*,4)$-face $f$ and
adjacent to two $3$-vertices not on $b(f)$.

\noindent{\rm (B5)}\ There is no $6$-vertex incident to two $(6,4^-,4^-)$-faces and one $(6,*,4)$-face.

\end{lemma}

\proof Let $L$ be a list assignment such that $|L(v)|=3$ for all $v\in V(G)$.
We make use of contradiction to show (B1)-(B5).
\begin{itemize}
\item [{\rm (B1)}] Suppose  that $v$ is adjacent to three $3$-vertices
$v_1, v_2$ and $v_3$.
Denote $G'=G-\{v,v_1,v_2,v_3\}$.
By the minimality of $G$,
$G'$ admits an $(L,1)^*$-coloring $\pi$.
Let $L_{\pi}$ be an induced list assignment of $G-G'$.
It is easy to deduce that $|L_{\pi}(v)|\geq 2$ and $|L_{\pi}(v_i)|\ge 1$ for each $i\in \{1,2,3\}$.
So for each $v_i$, we assign the color $\pi(v_i)\in L_{\pi}(v_i)$ to it.
Now we observe that there exists a color in $L_{\pi}(v)$ appearing at most once on the set $\{v_1, v_2, v_3\}$.
We color $v$ with such a color.
The obtained coloring is an $(L,1)^*$-coloring of $G$.
This contradicts the choice of $G$.

\item [{\rm (B2)}]  It suffices to prove that $G$ does not contain a $(4,4,4)$-face by (A3).
Suppose $f=[v_1v_2v_3]$ is a $3$-face with $d(v_1)=d(v_2)=d(v_3)=4$.
For each $i\in \{1,2,3\}$,
let $x_i, y_i$ denote the other two neighbors of $v_i$ not on $b(f)$.
Denote by $G'$ the graph obtained from $G$ by deleting edge $v_1v_2$.
By the minimality of $G$,
$G'$ has an $(L,1)^*$-coloring $\pi$.
If $\pi(v_1)\neq \pi(v_2)$,
then $G$ itself is $(L,1)^*$-colorable and thus we are done.
Otherwise, suppose $\pi(v_1)=\pi(v_2)$.
If $\pi$ is not an $(L,1)^*$-coloring of the whole graph $G$,
then without loss of generality,
assume that $\pi(v_1)=\pi(v_2)=\pi(x_1)=1$ and $\pi(v_3)=2$.
Moreover, none of $x_1$'s neighbors except $v_1$ is colored with 1.
First, we recolor each $v_i$ with a color $\pi'(v_i)$ in
$L(v_i)\setminus \{\pi(x_i), \pi(y_i)\}$, where $i\in \{1,2,3\}$.
We should point out that $\pi'(v_i)$ may be the same as $\pi(v_i)$,
but it does not matter.
Note that if at most two of $\pi'(v_1), \pi'(v_2), \pi'(v_3)$ are equal
then the resulting coloring is an $(L,1)^*$-coloring and thus we are done.
Otherwise, suppose that $\pi'(v_1)=\pi'(v_2)=\pi'(v_3)$.
Since $\pi'(v_1)\neq 1$ and $1\in L(v_1)$, we may further reassign color 1 to $v_1$ to
obtain an $(L, 1)^*$-coloring of $G$. This contradicts the choice of $G$.

\item [{\rm (B3)}]  Suppose $f=[v_1v_2v_3]$ is a $(5^+,4,4)$-face incident
to two light $4$-vertices $v_2$ and $v_3$.
By definition, we see that each $v_i$ ($i\in \{2,3\}$) is incident to two other $3$-vertices,
denoted by $x_i$ and $y_i$, which are not on $b(f)$.
Let $G'$ denote the graph obtained from $G$ by deleting edge $v_2v_3$.
Obviously, $G'$ has an $(L,1)^*$-coloring $\pi$ by the minimality of $G$.
Similarly,
if $\pi(v_2)\neq \pi(v_3)$,
then $G$ itself is $(L,1)^*$-colorable and thus we are done.
Otherwise, suppose $\pi(v_2)=\pi(v_3)$.
If $\pi$ is not an $(L,1)^*$-coloring of $G$,
then w.l.o.g.,
assume that $\pi(v_2)=\pi(v_3)=\pi(x_2)=1$ and $\pi(v_1)=2$.
Erase the color of $v_2$ and recolor $y_2$ with a color $a\in L(y_2)$ different from its neighbors.
If $L(v_2)\neq \{1,2,a\}$, then color $v_2$ with a color in $L(v_2)\setminus \{1,2,a\}$.
Otherwise, color $v_2$ with $a$. It is easy to verify that the
resulting coloring is an $(L,1)^*$-coloring of $G$, which is a contradiction.

\item [{\rm (B4)}]  Suppose  that a $5$-vertex $v$ is
incident to a $(5,*,4)$-face $f_1=[vv_1v_2]$ and adjacent to
two $3$-vertices $v_3$ and $v_4$.
Let $G'=G-\{v, v_1, v_2, v_3, v_4\}$.
By the minimality of $G$,
$G'$ has an $(L,1)^*$-coloring $\pi$.
Let $L_{\pi}$ be an induced list assignment of $G-G'$.
Obviously, $|L_{\pi}(v_i)|\ge 1$ for each $i\in \{1,\cdots,4\}$ and $|L_{\pi}(v)|\ge 2$.
By Lemma \ref{key-0}, $\pi$ can be extended to $G$, which is a contradiction.

\item [{\rm (B5)}]  Suppose that a $6$-vertex $v$ is incident to two $(6, 4^-, 4^-)$-faces $f_1, f_3$
and one $(6,*,4)$-face $f_5$ such that $d(v_i)\le 4$ for each $i=\{1,2,3,4\}$,
$d(v_6)=4$ and $v_5$
is an $\mathcal{S}$-vertex. Namely, $v_5$ is either a $3$-vertex or a light $4$-vertex.
Let $G'=G-\{v, v_1, v_2,\cdots, v_6\}$.
By minimality, $G'$ admits an $(L,1)^*$-coloring $\pi$.
Denote by $L_{\pi}$ an induced list assignment of $G-G'$.
It is easy to verify that
$|L_{\pi}(v_i)|\ge 1$ for each $i\in \{1,\cdots,6\}$ and $|L_{\pi}(v)|\ge 3$.
So we can color $v_i$ with $\pi(v_i)\in L_{\pi}(v_i)$ for each $i\in \{1,2,\cdots, 6\}$.
If there exists a color $a\in L_{\pi}(v)$ appearing at most once on the set $\{v_1, v_2, \cdots, v_6\}$,
then we further  assign color $a$ to $v$ and thus obtain an $(L,1)^*$-coloring of $G$.
Otherwise, each color in $L_{\pi}(v)$ appears exactly twice on the set $\{v_1, v_2, \cdots, v_6\}$.
Since $v_5$ is an $\mathcal{S}$-vertex, we can apply versions of arguments (i) and (ii) in the proof of Lemma \ref{key-0} to
obtain an $(L,1)^*$-coloring of $G$. \qed

\end{itemize}

\begin{lemma}\label{4-face-1}
Suppose that $f=[uvxy]$ is a $(3,4,m,4)$-face. Then

\noindent {\rm{(F1)}} $m\neq 3$.

\noindent {\rm{(F2)}} $x$ cannot be a soft $4$-vertex.

\end{lemma}
\proof (F1)\ Suppose to the contrary that $m=3$.
Let $G'=G-\{u,v,x,y\}$. By the minimality of $G$, $G'$ admits an $(L,1)^*$-coloring $\pi$.
Let $L_{\pi}$ be an induced list assignment of $G-G'$.
Notice that $|L_{\pi}(y)|\geq 1$, $|L_{\pi}(v)|\geq 1$,
$|L_{\pi}(u)|\geq 2$ and $|L_{\pi}(x)|\geq 2$.
First, we color $v$ with $a\in L_{\pi}(v)$ and
color $y$ with $b\in L_{\pi}(y)$.
Then color $u$ with $c\in L_{\pi}(u)\setminus \{a\}$ and
$x$ with $d\in L_{\pi}(x)\setminus \{b\}$.
One can easily check that the resulting coloring of $G$ is an $(L,1)^*$-coloring.
This contradicts the assumption of $G$.

(F2)\ Suppose to the contrary that $x$ is a soft $4$-vertex. By definition,
$x$ has other two neighbors whose degree are both 3, say $x_1$ and $x_2$.
Observe that neither $x_1$ nor $x_2$ is on $b(f)$.
Let $G'=G-\{u,v,x,y,x_1,x_2\}$.
Obviously, $G'$ admits an $(L,1)^*$-coloring $\pi$.
Let $L_{\pi}$ be an induced list assignment of $G-G'$.
For each $w\in \{v, y, x_1, x_2\}$, we deduce that $|L_{\pi}(w)|\ge 1$.
Moreover, $|L_{\pi}(u)|\ge 2$.
We first color $w$ with $\pi(w)\in L_{\pi}(w)$ and color $u$ with a color in $L_{\pi}(u)\setminus \{\pi(v)\}$.
If at least one of $x_1$ and $x_2$ has the same color as $\pi(v)$,
we can color $x$ with a color different from that of $v$ and $y$.
Otherwise, we can color $x$ with a color different from $x_1$ and $y$.
Therefore, we achieve an $(L,1)^*$-coloring of $G$, which is a contradiction. \qed

\begin{figure}[ht]
\centering
\includegraphics[width=2.1in]{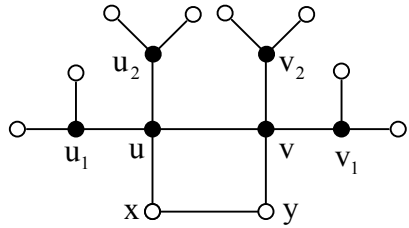}\\
\caption{Adjacent soft $4$-vertices $u$ and $v$.}
\label{fig2-1}
\end{figure}

\begin{lemma}\label{adjacent soft}
There is no adjacent soft $4$-vertices.
\end{lemma}

\proof Suppose to the contrary that $u$ and $v$ are adjacent soft $4$-vertices
such that $[uxyv]$ is a $4$-face and $u_1,u_2, v_1,v_2$ are $3$-vertices,
which is depicted in Figure \ref{fig2-1}.
By Observation \ref{three structures}(b), $u_i$ cannot be coincided with $v_j$, where $i, j\in \{1,2\}$.
Let $G'=G-\{u_1, u_2, v_1, v_2,u,v\}$.
For each $i\in \{1,2\}$, we color $u_i$ and $v_i$ with a color in $L_{\pi}(u_i)$ and $L_{\pi}(v_i)$, respectively.
If $L(u)\neq \{\pi(x), \pi(u_1), \pi(u_2)\}$,
then color $u$ with $a\in L(u)\setminus \{\pi(x), \pi(u_1), \pi(u_2)\}$.
It is easy to see that there exists at least one color in $L(v)\setminus \{\pi(y)\}$
which appears at most once on the set $\{u, v_1, v_2\}$. So we may assign such a color to $v$.
Now suppose that $L(u)=\{\pi(x), \pi(u_1), \pi(u_2)\}$.
By symmetry, we may suppose that $L(v)=\{\pi(y), \pi(v_1), \pi(v_2)\}$.
This implies that $\pi(v_1)\neq \pi(v_2)$. Thus, we can first color $u$ with $\pi(u_1)$ and
then assign a color in $L(v)\setminus \{\pi(u_1), \pi(y)\}$ to $v$.\qed

\begin{lemma}\label{lem-2}
Suppose $v$ is a $5$-vertex incident to two $3$-faces $f_1=[vv_1v_2]$
and $f_3=[vv_3v_4]$. Let $v_5$ be the neighbour of $v$ not belonging
to $f_1$ and $f_3$.
Then the following holds.

\noindent{\rm (C1)} If $f_1$ and $f_3$ are both $(5, 4^-, 4^-)$-faces, then $d(v_5)\ge 4$.

\noindent{\rm (C2)} If $f_1$ is a $(5,*,4)$-face and $f_3$ is a $(5,*,4^+)$-face, then $d(v_5)\ge 4$.

\noindent{\rm (C3)} $f_1$ and $f_3$ cannot be both $(5,*,4)$-faces.

\end{lemma}

\proof In each of following cases, we will show that an $(L,1)^*$-coloring of $G'\subset G$
can be extended to $G$, which is a contradiction.
\begin{itemize}
\item [{\rm (C1)}]  We only need to show that $d(v_5)\neq 3$ since $\delta(G)\ge 3$ by (A1).
Suppose that $v_5$ is a $3$-vertex.
Let $G'=G-\{v,v_1,\cdots,v_5\}$.
By the minimality of $G$,
$G'$ has an $(L,1)^*$-coloring $\pi$.
Let $L_{\pi}$ be an induced list assignment of $G-G'$.
It is easy to deduce that $|L_{\pi}(v_i)|\ge 1$ for each $i\in \{1,\cdots,5\}$ and $|L_{\pi}(v)|\ge 3$.
So we first color each $v_i$ with $\pi(v_i)\in L_{\pi}(v_i)$.
Observe that there exists a color $a\in L_{\pi}(v)$  that appears at most once on the set $\{v_1, v_2, \cdots, v_5\}$.
Therefore, we can color $v$ with $a$ to obtain an $(L,1)^*$-coloring of $G$.

\item [{\rm (C2)}]  Suppose that $d(v_2)=4$, $d(v_5)=3$ and $v_1$ and $v_3$ are both
$\mathcal{S}$-vertices. By definition, we see that $v_i$ is either a $3$-vertex
or a light $4$-vertex, where $i\in \{1,3\}$.
Let $G'=G-\{v, v_1, v_2, v_3, v_5\}$.
By the minimality of $G$,
$G'$ has an $(L,1)^*$-coloring $\pi$.
Let $L_{\pi}$ be an induced list assignment of $G-G'$.
The proof is split into two cases in light of the conditions of $v_3$.
\begin{itemize}
\item Assume $v_3$ is a $3$-vertex.
It is easy to calculate that $|L_{\pi}(v_i)|\ge 1$ for each $i\in \{1,2,3,5\}$ and $|L_{\pi}(v)|\ge 2$.
By Lemma \ref{key-0}, $\pi$ can be extended to $G$.

\item Assume $v_3$ is a light $4$-vertex.
By definition, let $x_3, y_3$ denote the other two neighbors of $v_3$
not on $b(f_3)$. Recolor $x_3$ and $y_3$ with a color different from its neighbors.
Next, we will show how to extend the resulting coloring $\pi'$ to $G$.
Denote $L_{\pi'}$ be the induced assignment of $G-G'$.
Notice that $|L_{\pi'}(v_i)|\ge 1$ for each $i\in \{1,2,5\}$.
If $|L_{\pi'}(v_3)|\ge 1$, then by Lemma \ref{key-0}, $\pi'$ can be extended to $G$.
Otherwise, we derive that $L(v_3)=\{\pi'(x_3), \pi'(y_3), \pi'(v_4)\}$.
First we assign a color in $L_{\pi'}(v_i)$ to each $v_i$, where $i\in \{1,2,5\}$.
It is easy to see that there is at least one color, say $a$,
belonging to $L(v)\setminus \{\pi'(v_4)\}$ that appears at most once on the set $\{v_1, v_2, v_5\}$.
We assign such a color $a$ to $v$. Then color $v_3$ with a color in $\{\pi'(x_3), \pi'(y_3)\}$ but different from $a$.
\end{itemize}

\item [{\rm (C3)}]  Suppose that $f_1$ and $f_3$ are both $(5,*,4)$-faces such
that $d(v_2)=d(v_4)=4$ and $v_1$ and $v_3$ are $\mathcal{S}$-vertices.
Let $G'=G-\{v, v_1, \cdots, v_4\}$.
Obviously, $G'$ has an $(L,1)^*$-coloring $\pi$ by the minimality of $G$.
Let $L_{\pi}$ be an induced list assignment of $G-G'$.
We assert that $v_i$ satisfies that $|L_{\pi}(v_i)|\ge 1$ for each $i\in \{1,\cdots, 4\}$
and $|L_{\pi}(v)|\ge 2$.
By Lemma \ref{key-0}, we can extend $\pi$ to the whole graph $G$ successfully. \qed

\end{itemize}

\begin{figure}[ht]
\centering
\includegraphics[width=2.1in]{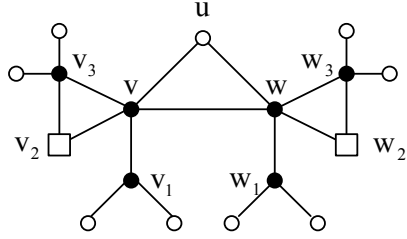}\\
\caption{The configuration  in Lemma \ref{key}.}
\label{fig4}
\end{figure}

\begin{lemma}\label{key}
There is no $3$-face incident to two bad $5$-vertices.
\end{lemma}
\proof  Suppose to the contrary that there is a $3$-face $[uvw]$ incident to two bad $5$-vertices $v$ and $w$,
depicted in Figure \ref{fig4}.
Let $G'=G-\{v,w,v_1,v_2,v_3,w_1,w_2,w_3\}$.
By the minimality of $G$,
$G'$ has an $(L,1)^*$-coloring $\pi$.
Let $L_{\pi}$ be an induced list assignment of $G-G'$.
Since each $w_i$ has at most two neighbors in $G'$,
we deduce that $|L_{\pi}(w_i)|\ge 1$ for each $i\in \{1,2,3\}$.
So we first color each $w_i$ with a color $\pi(w_i)\in L_{\pi}(w_i)$.
If $|L_{\pi}(w)|\ge 1$, namely $L(w)\neq \{\pi(u), \pi(w_1), \pi(w_2), \pi(w_3)\}$,
then by Lemma \ref{key-0} we may easy extend $\pi$ to $G$, since
$|L_{\pi}(v_i)|\ge 1$ for each $i\in \{1,2,3\}$.
Otherwise, we deduce that there exists a color $a$ in $L(w)\setminus \{\pi(u)\}$
that is the same as $\pi(w_{i^*})$ for some fixed $i^*\in \{1,2,3\}$. Color $w$ with $a$ and 
$v_i$ with a color $\pi(v_i)\in L_{\pi}(v_i)$ firstly, where $i\in \{1,2,3\}$.
For our simplicity, denote $V^*=\{v_1, v_2, v_3, w\}$.

First, suppose that there is a color, say $b\in L(v)\setminus \{\pi(u)\}$, appearing
at most once on the set $V^*$. We assign such a color $b$ to $v$.
If $b\neq a$, the obtained coloring is obvious an $(L,1)^*$-coloring.
Otherwise, assume that $b=a$. Now we erase the color $a$ from $w$.
One may check that the resulting coloring, say $\pi'$, satisfies that each of $v, w_1, w_2, w_3$ 
has at least one possible color in $G-G'$. In other words, $|L_{\pi'}(s)|\ge 1$ for each
$s\in \{v, w_1, w_2, w_3\}$. Hence, by Lemma \ref{key-0}, we can easily extend $\pi'$ to $G$.

Now, w.l.o.g., suppose that $L(v)=\{1,2,3\}$,
$\pi(u)=1$, $\pi(w)=2$ and each color in $\{2,3\}$ appears exactly twice on the set $V^*$.
It implies that $\pi(v_1)\in \{2,3\}$.
We apply versions of discussion (i) and (ii) in the proof of Lemma \ref{key-0}.
After doing that, one may check that now $v$ is colored with $\pi(v_2)$ and  
$v_1$ is recolored with a new color, say $\alpha$.
There are two cases left to discuss: if $\pi(v_2)=3$, namely the new color of $v$ is $3$, 
then the obtained coloring is an $(L,1)^*$-coloring and thus we are done;
otherwise, we uncolor $w$. Again, it is easy  
to see that the resulting coloring, say $\pi''$, satisfies that 
$|L_{\pi''}(s)|\ge 1$ for each
$s\in \{v, w_1, w_2, w_3\}$. 
Therefore, we can easily extend $\pi''$ to $G$ successfully by Lemma \ref{key-0}. \qed

\subsection{Discharging progress}

We now apply a discharging procedure to reach a contradiction.
Suppose that $u$ is adjacent to a $3$-vertex $v$ such that $uv$ is not incident to
any $3$-faces. We call $v$ a {\em free} $3$-vertex if $t(v)=0$
and a {\em pendant} $3$-vertex if $t(v)=1$.
For simplicity, we use $\nu_3(u)$ to denote the number of
free $3$-vertices adjacent to $u$ and $p_3(u)$ to denote the number of
pendant $3$-vertices of $u$. 
Suppose that $v$ is a soft $4$-vertex such that $f_1=[vv_1uv_2]$ is a $4$-face and 
$d(v_3)=d(v_4)=3$. If the opposite face to $f_1$ via $v$, i.e., $f_3$, is of degree at least $5$, 
then we call $v$ a {\em weak} $4$-vertex. We notice that 
every weak $4$-vertex is soft but not vice versa.

For $x\in V(G)$ and $y\in F(G)$, let $\tau(x\to y)$ denote the amount of
weights transferred from $x$ to $y$.
Suppose that $f=[v_1v_2v_3]$ is a 3-face. We use
$(d(v_1),d(v_2),d(v_3))\to (c_1,c_2,c_3)$ to denote $\tau(v_i\to f)=c_i$ for $i=1,2,3$.
Our discharging rules are defined as follows:

\noindent{\bf\rm {(R1)}} Let $f=[v_1v_2v_3]$ be a $3$-face.  We set

\ \ \ \ \noindent{\bf\rm {(R1.1)}} $(3,4,5^+)\to (0, 1, 3);$

\ \ \ \ \noindent{\bf\rm {(R1.2)}} $(3,5^+,5^+)\to (0, 2, 2);$

\ \ \ \ \noindent{\bf\rm {(R1.3)}}

\noindent
\[
(4,4,5^+)\rightarrow \left\{
\begin{array}{ll}
(0,1,3) & \mbox{if $v_1$ is a light $4$-vertex; \ \ \ \ \ \ \ \ \ \ \ \ \ \ \ \ \ \ \ \ \ }\medskip\\
(1,1,2) & \mbox{if neither $v_1$ nor $v_2$ is a light $4$-vertex.}
\end{array}\right.
\]

\ \ \ \ \noindent{\bf\rm {(R1.4)}}

\noindent
\[
(4,5^+,5^+)\rightarrow \left\{
\begin{array}{ll}
(1,1,2) & \mbox{if $v_2$ is a bad $5$-vertex; \ \ \ \ \ \ \ \ \ \ \ \ \ \ \ \ \ \ \ \ \ }\medskip\\
(0,2,2) & \mbox{if neither $v_2$ nor $v_3$ is a bad $5$-vertex.}
\end{array}\right.
\]

\ \ \ \ \noindent{\bf\rm {(R1.5)}}
\[
(5^+,5^+,5^+)\rightarrow \left\{
\begin{array}{ll}
(1,\frac 3 2,\frac 3 2) & \mbox{if $v_1$ is a bad $5$-vertex; \ \ \ \ \ \ \ \ \ \ \ \ \ \ \ \ \ \ \ \ \ }\medskip\\
(\frac 4 3, \frac 4 3, \frac 4 3) &   \mbox{if none of $v_1, v_2, v_3$ is a bad $5$-vertex.}
\end{array}\right.
\]

\noindent{\bf\rm {(R2)}} Suppose that $v$ is a $5^+$-vertex incident to a $4$-face $f=[vv_1uv_2]$. Then

\ \ \ \ \noindent{\bf\rm {(R2.1)}} $\tau(v\to f)=1$ if $d(v_1)\ge 4$ and $d(v_2)\ge 4$;

\ \ \ \ \noindent{\bf\rm {(R2.2)}} $\tau(v\to f)=\frac 4 3$ otherwise.

\noindent{\bf\rm {(R3)}} Suppose that $v$ is a non-weak $4$-vertex incident to a $4$-face $f=[vv_1uv_2]$.

\ \ \ \ \noindent{\bf\rm {(R3.1)}} Assume $d(v_1)=d(v_2)=3$. Then

\ \ \ \ \ \ \ \ \noindent{\bf\rm {(R3.1.1)}}  $\tau(v\to f)=\frac 4 3$ if the opposite face to $f$ via $v$ is of degree $3$;

\ \ \ \ \ \ \ \ \noindent{\bf\rm {(R3.1.2)}}  $\tau(v\to f)=\frac 2 3$ otherwise.

\ \ \ \ \noindent{\bf\rm {(R3.2)}}  Assume $d(v_1)\ge 4$ and $d(v_2)\ge 4$. Then

\ \ \ \ \ \ \ \ \noindent{\bf\rm {(R3.2.1)}}  $\tau(v\to f)=1$ if at least one of $v_1$ and $v_2$ is a soft $4$-vertex;

\ \ \ \ \ \ \ \ \noindent{\bf\rm {(R3.2.2)}}  $\tau(v\to f)=\frac 2 3$ otherwise.

\ \ \ \ \noindent{\bf\rm {(R3.3)}} Assume $d(v_1)=3$ and $d(v_2)\ge 4$. Then $\tau(v\to f)=\frac 2 3$.

\noindent{\bf\rm {(R4)}} Every $4^+$-vertex sends $1$ to each pendant $3$-vertex and
$\frac 1 3$ to each free $3$-vertex.

\medskip

According to (R3), we notice that a weak $4$-vertex does not send any charge.\\
We first consider the faces. Let $f$ be a $k$-face.

{\bf {Case}} $k=3$. Initially $\ch(f)=-4$.
Let $f=[v_1v_2v_3]$ with $d(v_1)\le d(v_2)\le d(v_3)$.
By (A1), $d(v_1)\ge 3$. If $d(v_1)=3$, then $d(v_2)\ge 4$ by (A2).
Together with (B2), we deduce that $f$ is either a $(3,4,5^+)$-face,
a $(3, 5^+, 5^+)$-face, a $(4, 4, 5^+)$-face, a $(4,5^+,5^+)$-face or a $(5^+, 5^+, 5^+)$-face.
It follows from (B3) and Lemma \ref{key} that
every possibility is indeed covered by rule (R1).
Obviously, $f$ takes charge 4 in total from its incident vertices.
Therefore,  $\nch(f)=-4+4=0$.

{\bf {Case}} $k=4$. Clearly, $w(f)=-2$.
Assume that $f=[vxuy]$ is a $4$-face.
By (A2), there are no adjacent $3$-vertices in $G$.
It follows that $f$ is incident to at most two $3$-vertices.
By symmetry, we have to discuss three cases depending on the conditions of these $3$-vertices.
\begin{itemize}
\item  $d(x)=d(y)=3$. By (F1), we deduce that
at least one of $u$ and $v$ is of degree at least $5$.
Moreover, if one of $u$ and $v$ is a $4$-vertex, say $v$, we claim that
$v$ cannot be weak by definition and (B1).
Hence,
$\nch(f)\ge -2+\frac 4 3+\frac 2 3=0$ by (R2) and (R3).

\item $d(x)=3$ and $d(y)\ge 4$. Note that $u$ and $v$ are both $4^+$-vertices.
Similarly, neither $u$ nor $v$ can be a weak $4$-vertex.
It follows from (R3.3) and (R2) that each of $u$ and $v$
sends charge at least $\frac 2 3$ to $f$.
So if one of them is a $5^+$-vertex, say $v$, then by (R2) we have that $\tau(v\to f)=\frac 4 3$ and
thus $f$ gets $\frac 2 3+\frac 4 3=2$ in total from incident vertices of $f$.
Otherwise, suppose $d(u)=d(v)=4$.
Now by (F2), $y$ cannot be a soft $4$-vertex and thus not weak. Hence,
$\nch(f)\ge -2+\frac 2 3\times 3=0$ by (R3.2).

\item $d(x)\ge 4$ and $d(y)\ge 4$.
Namely, $f$ is a $(4^+, 4^+, 4^+, 4^+)$-face.
If at most one of $u,v,x,y$ is a weak $4$-vertex, then
$\nch(f)\ge -2+\frac 2 3\times 3=0$.
Otherwise, by Lemma \ref{adjacent soft}, assume that $v$ and $u$ are weak $4$-vertices and thus soft.
We see that $\tau(x\to f)=\tau(y\to f)=1$ by (R3.2.1) and (R2.1) which implies that $\nch(f)\ge -2+1\times 2=0$.

\end{itemize}

{\bf {Case}} $k\ge 5$. Then $\nch(f)=\ch(f)=2d(f)-10\ge 0$.

\medskip
Now we consider the vertices.
Let $v$ be a $k$-vertex with $k\ge 3$ by (A1).
For $v\in V(G)$, we use $m_4(v)$ to denote the number of $4$-faces incident to $v$.
So by Observation \ref{three structures} (a) and (b),
we derive that
$t(v)\leq \lfloor \frac {d(v)}{2}\rfloor$ and $m_4(v)\leq \lfloor \frac {d(v)}{2}\rfloor$.
Furthermore, $t(v)+m_4(v)\le \lfloor \frac {d(v)}{2}\rfloor$ by Observation \ref{three structures} (c).

\begin{observation}\label{charge}
Suppose $v$ is a $4^+$-vertex which is incident to a $3$-face $f$. Then, by (R1), we have the following:

\noindent{\rm {(a)}} $\tau(v\to f)\le 1$ if $d(v)=4$;

\noindent{\rm {(b)}} $\tau(v\to f)\in \{3, 2, \frac 3 2, \frac 4 3, 1\}$ if $d(v)\ge 5$;
moreover, if $\tau(v\to f)=3$ then $f$ is a $(5^+,*,4)$-face.

\end{observation}

{\bf {Case}} $k=3$. Then $\ch(v)=-1$.
Clearly, $t(v)\le 1$. If $t(v)=1$, then there exists a neighbor of $v$, say $u$, so that
$v$ is a pendant $3$-vertex of $u$. By (A2), $d(u)\ge 4$.
Thus, $\nch(v)=-1+1=0$ by (R4).
Otherwise, we obtain that $\nch(v)=-1+\frac 1 3\times 3=0$ by (R4).

{\bf {Case}} $k=4$. Then $\ch(v)=2$.
Note that $t(v)\le 2$.
If $t(v)=2$, then $m_4(v)=0$ and $p_3(v)=0$.
So $\nch(v)\ge 2-1\times 2=0$ by Observation \ref{charge} (a).
If $t(v)=0$, then $n_3(v)\le 2$ by (B1) and $m_4(v)\le 2$. We need to consider following cases.
\begin{itemize}
\item  $m_4(v)=2$. W.l.o.g., assume that $f_1=[vv_1uv_2]$ and $f_3=[vv_3wv_4]$ are incident $4$-faces.
Obviously, $p_3(v)=0$ by Observation \ref{three structures} (b).
However, $\nu_3(v)\le 2$ by (B1).
By (R3), $v$ sends charge at most 1 to $f_i$, where $i=1,3$.
If $n_3(v)=0$, then $\nu_3(v)=0$ and thus $\nch(v)\ge 2-1\times 2=0$. If $n_3(v)=1$, say $v_1$ is a $3$-vertex,
then $\tau(v\to f_1)\le \frac 2 3$ by (R3.3) and thus $\nch(v)\ge 2-\frac 2 3-1-\frac 1 3=0$ by (R4).
Now suppose that $n_3(v)=2$. By symmetry, we have two cases depending on the conditions of
these two $3$-vertices. If $d(v_1)=d(v_2)=3$, then $\tau(v\to f_1)=\frac 2 3$ by (R3.1.2).
By (B1), $v_3$ and $v_4$ are both $4^+$-vertices. 
Moreover, neither $v_3$ nor $v_4$ is a soft $4$-vertex according to Lemma \ref{adjacent soft}. 
So by (R3.2.2), $\tau(v\to f_3)\le \frac 2 3$.
Hence $\nch(v)\ge 2-\frac 2 3-\frac 2 3-\frac 1 3\times 2=0$.
Otherwise, suppose that $d(v_i)=d(v_j)=3$, where $i\in \{1,2\}$ and $j\in \{3,4\}$. 
We derive that 
$\nch(v)\ge 2-\frac 2 3\times 2-\frac 1 3\times 2=0$ by (R3.3).

\item  $m_4(v)=1$. W.l.o.g, assume that $d(f_1)=4$.
This implies that $d(f_3)\ge 5$.
Again, $\tau(v\to f_1)\le 1$ by (R3).
If $n_3(v)\le 1$ then we have that $\nch(v)\ge 2-1-1=0$ by (R4).
So in what follows, we assume that $n_3(v)=2$. If $d(v_3)=d(v_4)=3$ then $v$ is a weak $4$-vertex,
implying that $v$ sends nothing to $f_1$. So $\nch(v)\ge 2-1\times 2=0$ by (R4).
If $d(v_1)=d(v_2)=3$, then $p_3(v)=0$ by Observation \ref{three
  structures} (b).
We deduce that $\nch(v)\ge 2-\frac 2 3-\frac 1 3 \times 2=\frac 2 3$ by (R3.1.2) and (R4).
Otherwise, suppose $d(v_i)=d(v_j)=3$, 
where $i\in \{1,2\}$ and $j\in \{3,4\}$. 
It follows immediately from (R3.3)and (R4) that
$\nch(v)\ge 2-\frac 2 3-1-\frac 1 3=0$.

\item $m_4(v)=0$. Obviously, $\nch(v)\ge 2-1\times 2=0$ by (R4).

\end{itemize}

Now, in the following,
we consider the case $t(v)=1$.
Assume that $f_1$ is a $3$-face.
By (A1) and (B2), $f_1$ is either a $(4,3,5^+)$-face, a $(4,4,5^+)$-face or
a $(4, 5^+, 5^+)$-face. Observe that $m_4(v)\le 1$. First assume that $m_4(v)=0$.
If $f_1$ is a $(4,3,5^+)$-face, then $p_3(v)\le 1$ by (B1) and
hence $\nch(v)\ge 2-1-1=0$ by Observation \ref{charge} (a) and (R2).
Next suppose that $f_1$ is a $(4,4,5^+)$-face.
If $n_3(v)=2$, then $v$ is a light $4$-vertex.
By (R1.3), we see that $v$ sends nothing to $f_1$ and therefore $\nch(v)\ge 2-1\times 2=0$ by (R4).
Otherwise, at most one of
$v_3, v_4$ is a $3$-vertex and hence $\nch(v)\ge 2-1-1=0$ by
Observation \ref{charge} (a) and (R4).
Finally, we suppose that $f_1$ is a $(4, 5^+, 5^+)$-face.
If neither $v_1$ nor $v_2$ is a bad $5$-vertex, then $v$
sends nothing to $f_1$ by (R1.4) and thus $\nch(v)\ge 2-1\times 2=0$ by (R4).
Otherwise, one of $v_1$ and $v_2$ is a bad $5$-vertex.
If follows directly from (C2) that $n_3(v)\le 1$.
Therefore, $\nch(v)\ge 2-1-1=0$ by (R2).
Now suppose that $m_4(v)=1$. By Observation \ref{three structures} (c),
we may assume that $f_3=[vv_3wv_4]$ is a $4$-face. In this case, $p_3(v)=0$.  
If $d(v_3)=d(v_4)=3$, then $\tau(v\to f_3)=\frac 4 3$ by (R3.1.1).
It follows from (B1) and (C2) that $f$ is neither a $(4,3,5^+)$-face nor a $(4,5,5^+)$-face such that $v_2$ is a bad $5$-vertex.
So we deduce that $f_1$ gets nothing from $v$ by (R1.3), which
implies that $\nch(v)\ge 2-\frac 4 3-\frac 1 3\times 2=0$.
If exactly one of $v_3, v_4$ is a $3$-vertex, then $\tau(v\to f_3)\le \frac 2 3$ by (R3,3).
Thus, $\nch(v)\ge 2-1-\frac 2 3-\frac 1 3=0$ by Observation \ref{charge} (a) and (R4).
Otherwise, we suppose that $v_3, v_4$ are both of degree at least 4.
In this case, $\nu_3(v)=0$ and hence $\nch(v)\ge 2-1-1=0$ by (R3.2) and Observation \ref{charge} (a).

{\bf {Case}} $k=5$. Then $\ch(v)=5$.
Also, $t(v)\le 2$. we have three cases to discuss.

Assume $t(v)=0$. If $m_4(v)=0$, then $\nch(v)\ge 5-1\times 5=0$ by (R4).
If $m_4(v)=1$, then  $p_3(v)\le 3$. Thus $\nch(v)\ge 5-\frac 4 3-1\times 3-2\times \frac 1 3=0$ by (R2) and (R4).
Now suppose that $m_4(v)=2$. By Observation \ref{three structures} (c), we assert that $p_3(v)\le 1$.
So $\nch(v)\ge 5-\frac 4 3\times 2-\frac 1 3\times 4-1=0$.

Next assume $t(v)=1$, say $f_1$. Then $\tau(v\to f_1)\le 3$ by Observation \ref{charge} (b).
Moreover, equality holds iff $f_1$ is a $(5,*,4)$-face. 
So if $\tau(v\to f_1)=3$ then at most one of $v_3, v_4, v_5$ is a $3$-vertex by (B4). 
Furthermore, $m_4(v)\le 1$.
When $m_4(v)=0$, 
we deduce that $\nch(v)\ge 5-3-1=1$ by (R4).
When $m_4(v)=1$, by symmetry, say $f_3$ is a $4$-face, we have two cases to discuss:
if $p_3(v)=1$, namely, $v_5$ is a $3$-vertex, then $\tau(v\to f_3)\le 1$ by (R2)
and neither $v_3$ nor $v_4$ takes charge from $v$. Thus $\nch(v)\ge 5-3-1-1=0$;
otherwise, $p_3(v)=0$ and we have $\nch(v)\ge 5-3-\frac 4 3-\frac 1 3=\frac 1 3$.
Now suppose that $\tau(v\to f_1)\le 2$.
By (R2) and (R4), $\nch(v)\ge 5-2-1\times 3=0$ if $m_4(v)=0$ and
$\nch(v)\ge 5-2-\frac 4 3-1-2\times \frac 1 3=0$ if $m_4(v)=1$.

Now assume $t(v)=2$. By symmetry, assume $f_1$ and
$f_3$ are both $3$-faces. Observe that $m_4(v)=0$.
For simplicity, denote $\tau(v\to f_1)=\sigma_1$ and $\tau(v\to f_3)=\sigma_2$.
Let $\sigma=\max \{\sigma_1, \sigma_2\}$.
If $\sigma\le 2$, then $\nch(v)\ge 5-2\times 2-1=0$ by (R2).
Now assume that $\sigma=3$, i.e., $f_1$ gets charge $3$ from $v$.
It means that $f_1$ is a $(5,*,4)$-face by Observation \ref{charge}.
By (C3), $f_3$ cannot be a $(5,*,4)$-face.
This implies that $\sigma_2\le 2$.
Moreover, if $v_5$ is a $3$-vertex, then $f_3$ is neither a $(5,*,4^+)$-face by (C2)
nor a $(5,4,4)$-face by (C1).
It follows from (R1.4) and (R1.5) that $\sigma_2\le 1$, since $v$ is a bad $5$-vertex.
Thus, $\nch(v)\ge 5-3-1-1=0$ by (R2).
Otherwise, we easily obtain that $\nch(v)\ge 5-3-2=0$.

{\bf {Case}} $k\ge 6$. Notice that  $t(v)\le \lfloor \frac {d(v)}{2}\rfloor$.
If $v$ is incident to a $4$-face $f_i$, then by (R2) we inspect
$v$ sends a charge at most $\frac 4 3$ to $f_i$, while $\frac 1 3$ to each of $v_i$ and $v_{i+1}$.
So we may consider $v$ as a vertex which sends charge at most $\frac 4 3+2\times \frac 1 3=2$ to $f_i$.
So by  (R4) and Observation \ref{charge}, we have
\begin{eqnarray*}
\nch(v)&\geq&  3d(v)-10-3t(v)-2m_4(v)-(d(v)-2t(v)-2m_4(v))\\
&=& 2d(v)-10-t(v)\equiv \tau(v)
\end{eqnarray*}

If $d(v)\ge 7$, then $\tau(v)\ge 2d(v)-10-\frac {d(v)}{2}=\frac 3 2 d(v)-10\ge \frac 3 2\times 7-10=\frac 1 2>0$.
Now suppose that $d(v)=6$. If $t(v)\le 2$ then $\tau(v)\ge 2\times 6-10-2=0$.
So, in what follows, assume that $t(v)=3$ and $d(f_i)=3$ for $i=1,3,5$.
Clearly, $m_4(v)=0$.
Similarly, if there are at most two of $3$-faces get charge $3\times 2$ in total from $v$,
then $\nch(v)\ge 8-2\times 3-2=0$.
Otherwise, suppose $\tau(v\to f_i)=3$ for each $i\in \{1,3,5\}$.
By Observation \ref{charge} (b), we assert that
$f_i$ is a $(6, *, 4)$-face.
Noting that a $(6, *, 4)$-face is also a $(6,4^-,4^-)$-face,
we may regard $v$ as a $6$-vertex which
is incident to two $(6,4^-,4^-)$-faces and one $(6, *, 4)$-face.
However, it is impossible by (B5).

Therefore, we complete the proof of Theorem \ref{main01}.\qed

\end{document}